\documentclass[11pt,leqno]{amsart}
\usepackage{amsmath}
\usepackage{amsthm}
\usepackage{amssymb}
\theoremstyle{definition}

\theoremstyle{plain}

\begin{document}
\setcounter{page}{177}
\begin{flushright}
\scriptsize{{\em Journal of the Nigerian Association of Mathematical Physics}\\
{\em Volume} {\bf 10}\,({\em November} 2006), 177--180 (Original Print)\\
\copyright{\em J of NAMP}
}
\end{flushright}
\vspace{16mm}

\begin{center}
{\normalsize\bf On the successive coefficients of certain univalent functions} \\[12mm]
     {\sc K. O. BABALOLA$^{1}$} \\ [8mm]

\begin{minipage}{100mm}
{\small {\sc Abstract}

$\frac{\rule{3.90in}{0.02in}}{{}}$\\
The object of this paper is to study relationship between successive coefficients of some subclasses of the class of univalent functions in the unit disk. the result obtained is sharp, and is used to provide a new, short proof of the well-known conjecture of Robertson on the coefficients of close-to-convex functions.
$\frac{{}}{\rule{3.90in}{0.02in}}$
}
\end{minipage}
\end{center}

 \renewcommand{\thefootnote}{}
 \footnotetext{$^1$\,Department of Mathematics, University of Ilorin, Ilorin, Nigeria. E-mail:
abuqudduus\symbol{64}yahoo.com}

\def\iff{if and only if }
\def\S{Smarandache }
\newcommand{\norm}[1]{\left\Vert#1\right\Vert}

\vskip 12mm

{\bf 1.0 Introduction}
\medskip

Let $A$ denote the class of functions of the form:
$$f(z)=z+a_2z^2+\cdots\eqno{(1.1)}$$
which are analyticin the unit disk $E=\{z\colon |z|<1\}$. Also let $P$ be the class of functions
$$p(z)=1+c_1z+c_2z^2+\cdots\eqno{(1.2)}$$
which are analytic in $E$ and have positive real part. It is well-known that for $f\in A$ the condition that
$$Re\frac{zf'(z)}{f(z)}>0,\;\;z\in E\eqno{(1.3)}$$
is necessary and sufficient for starlikeness (and univalence) in the unit disk. Also necessary and sufficient for $f\in A$ to be convex in the unit disk is that
$$Re\left(1+\frac{zf''(z)}{f'(z)}\right)>0.\eqno(1.4)$$

The families of functions, denoted respectively by $S^\ast$ and $C$, were discovered by Robertson \cite{4}. They have attracted attention of reasearchers in geometric functions theory, and have been generalized in \cite{6} by Salagean, who said a function $f\in A$ belongs to the class $S_n(\alpha)$, $0\leq\alpha<1$, if and only if
$$Re\frac{D^{n+1}f(z)}{D^nf(z)}>\alpha,\;\;z\in E\eqno(1.5)$$
where $n\in N_0=\{0,1,2,\cdots\}$, $0\leq\alpha<1$ and $D^n$ is defined as follows:
$$D^nf(z)=D(D^{n-1}f(z))=z[D^{n-1}f(z)]'$$
with $D^0f(z)=f(z)$.\vskip 2mm

It can be observed that the cases $n=0$ and $n=1$ respectively correspond to the families of starlike and convex functions of order $\alpha$ in $E$. The Salagean derivative has gained much acclaim as a unifying factor in the study of many classes of functions. For the case $\alpha=0$, we shall simply write $S_n$.\vskip 2mm

The purpose of the present studyis toobtain the best possible relationship between successive coefficients of functions of the subclasses, $S_n$. We express the extremal functions in terms of the integral operator, $I_n$, also introduced in \cite{6} as follows:
$$I_nf(z)=I(I_{n-1}f(z))=\int_0^z\frac{I_{n-1}f(t)}{t}dt$$
with $I_0f(z)=f(z)$.\vskip 2mm

In section 2 we state the basic results on which we shall depend for the proof of our result in section 3. Our proof follows a method devised by Leung \cite{2} for the proof of the case $n=0$. As a consequence of the main result we provide a new short proof of a conjecture of Robertson on the coefficients of close-to-convex univalent functions  \cite{5}.
 \medskip

{\bf 2.0 Preliminary Lemmas}\vskip 2mm

In section 3, we shall be making use of the following lemmas:\vskip 2mm

{\bf Lemma 2.1}(\cite{1})\vskip 2mm

{\em Let $\phi(z)=\sum_{j=0}\lambda_jz^j$ be an arbitrary power series having a positive radius of convergence and normalized by $\phi(0)=0$. Also let $\exp(\phi(z))=\sum_{j=0}\beta_jz^j$, $\beta_0=1$, be the exponentiated power series of $\phi(z)$. Then

$$|\beta_k|\leq\exp\left\{\sum_{j=1}^k\left(j|\lambda_j|^2-\frac{1}{j}\right)\right\}.$$}

{\bf Lemma 2.2}(\cite{2})\vskip 2mm

{\em For every $p\in P$ and every positive integer $k$, there exists a complex number $\nu$ with $|\nu|=1$ such that:
$$\sum_{j=1}^k\frac{1}{j}\left|p_j-\nu^j\right|^2\leq\sum_{j=1}^k\frac{1}{j}.$$}

\medskip

{\bf 3.0 Main Result}\vskip 2mm

{\bf Theorem 3.1}\vskip 2mm

{\em Let $f\in S_n$, then for any complex number $\nu$ such that $|\nu|=1$,
$$|(k+1)^na_{k+1}-\nu k^na_k|\leq 1;\;\;k=1,2,3,\cdots\eqno(3.1)$$
For any fixed $n\in N_0$, the inequality is the best possible, with equality only for the function:
$$f_n(z)=I_n\left\{\frac{z}{(1-\nu z)(1-\gamma z)}\right\},\;\;|\gamma|=1.\eqno(3.2)$$
}

By triangle inequality, the first corollary below is immediate and the next two are its consequences.\vskip 2mm

{\bf Corollary 3.2}\vskip 2mm

{\em Let $f\in S_n$, then
$$|(k+1)^n|a_{k+1}|-k^n|a_k||\leq 1;\;\;k=1,2,3,\cdots$$
}

{\bf Corollary 3.3}\vskip 2mm

{\em For every $f\in S_n$
$$|a_k|\leq k^{1-n};\;\;k=2,3,\cdots$$
}

{\bf Corollary 3.3}\vskip 2mm

{\em For every odd function $f\in S_n$
$$|a_{2k+1}|\leq (2k+1)^{-n};\;\;k=1,2,3,\cdots$$
}

{\bf Proof [Theorem 3.1]}\vskip 2mm

Take $\alpha=0$ in equation (1.5). Then for some $p\in P$, we have
$$\frac{D^{n+1}f(z)}{D^nf(z)}=p(z)\eqno(3.3)$$
which implies
$$\frac{[D^nf(z)]'}{D^nf(z)}=\frac{1}{z}+\sum_{j=1}p_jz^{j-1}\eqno(3.4)$$
Integrating (3.4) we obtain
$$\log_e\frac{D^nf(z)}{z}=\sum_{j=1}\frac{p_j}{j}z^j\eqno(3.5)$$
Now, for $|\nu|=1$, we have
$$\log_e\left\{(1-\nu z)\frac{D^nf(z)}{z}\right\}=\sum_{j=1}\frac{1}{j}(p_j-\nu^j)z^j\eqno(3.6)$$
whereas
$$(1-\nu z)\frac{D^nf(z)}{z}=\sum_{j=0}\{(1+j)^na_{j+1}-\nu j^na_j\}z^j.\eqno(3.7)$$
Combining equations (3.6) and (3.7) we get
$$\sum_{j=0}\{(1+j)^na_{j+1}-\nu j^na_j\}z^j=\exp\left\{\sum_{j=1}\frac{1}{j}(p_j-\nu^j)z^j\right\}.\eqno(3.8)$$
Appling Lemma 2.2 to (3.8) we have
$$|(k+1)^na_{k+1}-\nu k^na_k|\leq\exp\left\{\sum_{j=1}^k\frac{1}{j}|p_j-\nu^j|^2-\frac{1}{j}\right\}.\eqno(3.8)$$

By Lemma 2.2, our choice of $\nu$ ensures that the exponent on the right of equation (3.9) is nonpositive. Hence we have our result.\vskip 2mm

Observe from (3.9) that equality occurs in (3.1) if and only if for any complex number $\gamma$ ($|\gamma|=1$),
$$p_j-\nu^j=\gamma^j.\eqno(3.10)$$
Using (3.10) in (3.8) we obtain
$$\sum_{j=0}\{(j+1)^na_{j+1}-\nu j^na_j\}z^j=\exp\left\{\sum_{j=1}\frac{\gamma^j}{j}z^j\right\}.\eqno(3.11)$$
Equation (3.11) implies 
$$(j+1)^na_{j+1}-\nu j^na_j=\gamma^j,\;\;j=1,2,3\cdots\eqno(3.12)$$

It can be shown from (3.12) using simple inductive argument that the coefficients of the extremal functions satisfy
$$(k+1)^na_{k+1}=\frac{\nu^{k+1}-\gamma^{k+1}}{\nu-\gamma}.$$

Thus the extremal functions are those given by (3.2). This concludes the proof.\vskip 2mm

{\bf Remark 3.5.}\vskip 2mm

{\em The case $n=0$ yields results which are due to Leung \cite{2} while for $n=1$, we have equivalent results for functions which map the unit disk onto a convex domain.}

\medskip

In the sequel we give an important consequence of Theorem 3.1, which is a new, concise proof of the Robertson's conjecture on the coefficient of close-to-convex univalent function. Let
$$g(z)=z+b_2z^2+\cdots\eqno(3.15)$$
be a convex function. Then a function $f(z)$ in $A$ is said to be close-to-convex in $E$ if and only if
$$Re\frac{f'(z)}{g'(z)}>0,\;\;z\in E.\eqno(3.16)$$

We denote the family of close-to-convex functions by $K$. We remark that the Robertson's conjecture, which is the following, has been proved by Leung in \cite{3}, but not without sweat.\vskip 2mm

{\bf Theorem 3.1}\vskip 2mm

{\em For every $f\in K$ and positive integers k, m;

$$|k|a_k|-m|a_m||\leq |k^2-m^2|.\eqno(3.17)$$
}

{\bf Proof}\vskip 2mm

Since $f\in K$, there exists $p\in P$ and a $g\in C$ such that $f'(z)=g'(z)p(z)$, $z\in E$. Hence we have
$$\sum_{k=0}(k+1)a_{k+1}z^k=\sum_{k=0}c_kz^k\eqno(3.18)$$
where
$$c_k=\sum_{j=0}^kp_j(k+1-j)b_{k+1-j}\eqno(3.19)$$
is the coefficient of the product series of $g'(z)$ and $p(z)$ with $p_0=b_1=1$.\vskip 2mm

If we multiply (3.18) by $1-z$ and compare coefficients, we get
$$(k+1)a_{k+1}-ka_k=c_k-c_{k-1}.\eqno(3.20)$$
Thus using (3.19) in (3.20) and applying triangle inequality, we obtain
$$|(k+1)|a_{k+1}|-k|a_k||\leq|(k+1)b_{k+1}-kb_k|+\sum_{j=1}^{k-1}|p_j||(k+1-j)b_{k+1}-(k-j)b_k|+|p_k|.\eqno(3.21)$$

Thus by Theorem 3.1, for any complex number $\nu$ such that $|\nu|=1$, $|(k+1)b_{k+1}-\nu kb_k|\leq 1$, $k=1,2,\cdots$ for $g\in C$. Using this together with the well known Caratheodory inequality $|p_k|\leq 2$, $k=1,2,\cdots$, equation (3.21) gives
$$|(k+1)|a_{k+1}|-k|a_k||\leq|2k+1.\eqno(3.22)$$

The proof is concluded by observing that general case of the inequality (3.17) follows from (3.22) by simple induction.\vskip 2mm

Strict inequality holds for all $k$ and $m$, unless $f(z)$ is a rotation of the Koebe function $k(z)=\tfrac{z}{(1-z)^2}$.

\medskip

{\bf 4.0 Conclusion}\vskip 2mm

The main results of this paper provides the basic relationship between the coefficients of certain subclasses of univalent functions, expressed in the best possible inequality. The result generalized an earlier on for starlike functions by Leung \cite{2}. A very important consequence followed from the main result, that is we provide a new, short proof of the old conjecture of Robertson \cite{5}. Although a proof of the conjecture had appeared earlier in \cite{3}, yet the conciseness of our method of proof nakes a significant contribution to the body of knowledge in geometric functions theory.

\end{document}